\documentclass[11pt,reqno]{amsart}
\usepackage[a4paper,margin=1.1in]{geometry}
\usepackage{amsmath,amssymb,amsfonts,amsthm,mathtools}
\usepackage[hidelinks]{hyperref}
\usepackage{url}
\allowdisplaybreaks

\newtheorem{theorem}{Theorem}
\newtheorem{corollary}[theorem]{Corollary}
\newtheorem{conjecture}[theorem]{Conjecture}
\newtheorem{proposition}[theorem]{Proposition}
\newtheorem{lemma}[theorem]{Lemma}
\theoremstyle{remark}
\newtheorem{remark}[theorem]{Remark}

\numberwithin{equation}{section}
\title[Lower bounds in a conjecture of Ng]{Lower bounds for discrete moments in a conjecture of Ng}
\author{Andrew Pearce-Crump}
\address{School of Mathematics, University of Bristol, Fry Building, Woodland Road, Bristol, BS8 1UG, United Kingdom}
\email{andrew.pearce-crump@bristol.ac.uk}
\date{}

\begin{document}

\begin{abstract}
Assuming the Riemann Hypothesis and the simplicity of the non-trivial zeros of the Riemann zeta function $\zeta(s)$, we prove lower bounds for the discrete moments $\sum_{0<\gamma\le T}|\zeta(2\rho)/\zeta'(\rho)|^{2k}$ for every fixed rational $k\ge1$, and for every fixed real $k\ge1$. At $k=1$ we give a new proof of the explicit bound, half of the main term conjectured by Ng. 
\end{abstract}

\maketitle

\section{Introduction}\label{sect:intro}

Let $\zeta(s)$ denote the Riemann zeta function and write $\rho=\tfrac12+i\gamma$ for its non-trivial zeros. In this paper we address the two problems in item number~5 of the recent Comparative Prime Number Theory (CPNT) Problem List \cite{CPNT_ProblemList}. They were suggested by Nathan Ng, who in his thesis \cite{NgThesis} made the following conjecture.

\begin{conjecture}[Ng]\label{conj:ng}
Assume the Riemann Hypothesis and that the non-trivial zeros of $\zeta(s)$ are simple. Then
\begin{equation*}
\sum_{0 < \gamma \le T} \left| \frac{\zeta(2\rho)}{\zeta'(\rho)} \right|^2 \sim \frac{T}{2\pi}
\end{equation*}
as $T \to \infty$.
\end{conjecture}

The assumption of the simplicity of the zeros ensures that the sum is well-defined.

\begin{remark}\label{rem:hypotheses-ng}
In a private communication (2026), Ng noted that in stating Conjecture~\ref{conj:ng} he implicitly assumed the Riemann Hypothesis and the simplicity of the zeros. He added that the Riemann Hypothesis may be removed, and that he expects the conjecture to remain true assuming only that the zeros are simple. The simplicity assumption can in turn be removed from the statement by passing to the residue formulation
\[
\sum_{0 < \gamma \le T} \left| \operatorname*{Res}_{s=\rho} \frac{\zeta(2s)}{\zeta(s)} \right|^2 \sim \frac{T}{2\pi}.
\]
Here the sum is over distinct zeros. The residue is well-defined regardless of the multiplicity of a zero $\rho$ and, at a simple zero, equals $\zeta(2\rho)/\zeta'(\rho)$; its squared modulus is the summand in Conjecture~\ref{conj:ng}. Thus the two formulations agree under simplicity. Without simplicity, the residue formulation is a conjectural extension; rarity of multiple zeros alone would not establish that their contribution is negligible.
\end{remark}

For Conjecture~\ref{conj:ng} itself the lower bound is explicit.

\begin{theorem}\label{thm:NgCon}
Assume the Riemann Hypothesis and that the non-trivial zeros of $\zeta(s)$ are simple. Then, for any fixed $\varepsilon > 0$,
\begin{equation*}
\sum_{0 < \gamma \le T} \left| \frac{\zeta(2\rho)}{\zeta'(\rho)} \right|^2 \ge \left(\frac{1}{4\pi} - \varepsilon\right)T
\end{equation*}
for $T$ sufficiently large.
\end{theorem}

This lower bound is half the conjectured size, and the situation is closely parallel to that of the reciprocal moments. Ng based his conjecture on techniques from Gonek's MSRI talk \cite{Gonek}. Gonek had announced in \cite{Gonek1989}, and developed in that talk, the following conjecture.

\begin{conjecture}[Gonek]\label{conj:gonek}
Assume the Riemann Hypothesis and that the non-trivial zeros of $\zeta(s)$ are simple. Then
\begin{equation*}
\sum_{0 < \gamma \le T} \frac{1}{|\zeta'(\rho)|^2} \sim \frac{3}{\pi^3}T
\end{equation*}
as $T \to \infty$.
\end{conjecture}

Milinovich and Ng \cite{MilinovichNg} proved a lower bound for Gonek's conjecture which also differs from the conjectured size by a factor of two.

\begin{theorem}[Milinovich--Ng]\label{thm:MN}
Assume the Riemann Hypothesis and that the non-trivial zeros of $\zeta(s)$ are simple. Then, for any fixed $\varepsilon > 0$,
\begin{equation*}
\sum_{0 < \gamma \le T} \frac{1}{|\zeta'(\rho)|^2} \ge \left(\frac{3}{2\pi^3} - \varepsilon\right)T
\end{equation*}
for $T$ sufficiently large.
\end{theorem}

Conjecture~\ref{conj:ng} is the case $k=1$ of a more general conjecture of Ng \cite[Conjecture~1.11]{NgPNET}. Stated in the exponent $2k$ used throughout this paper, it predicts, for $0\le k<\tfrac32$,
\[
\sum_{0<\gamma\le T}\left|\frac{\zeta(2\rho)}{\zeta'(\rho)}\right|^{2k}
\sim\frac{G^2(2-k)}{G(3-2k)}\,a(-k)
\Bigl(\sum_{n=1}^{\infty}\frac{d_{k}(n)}{n^2}\Bigr)
\frac1{2\pi}\int_{2\pi}^T\Bigl(\log\frac t{2\pi}\Bigr)^{(k-1)^2}\,dt.
\]

Here $d_z(n)$ is the generalised divisor function and $G(z)$ is the Barnes $G$-function. The arithmetic factor is the one appearing in the conjecture of Hughes, Keating and O'Connell \cite{HKO} for the positive moments $\sum|\zeta'(\rho)|^{2k}$, namely
\[
 a(z)=\prod_p\left(1-\frac1p\right)^{z^2}
 \sum_{m=0}^{\infty}\left(\frac{\Gamma(m+z)}{m!\,\Gamma(z)}\right)^2p^{-m}.
\]
The range $0\le k<\tfrac32$ is inherited from the range $\Re(k)>-\tfrac32$ of that conjecture. Taking $k=1$ recovers Conjecture~\ref{conj:ng}.

Item number~5, Problem~2 of the CPNT Problem List asks for a lower bound of this form, namely
\begin{equation}\label{eq:Ng-higher-target}
  \sum_{0<\gamma\le T}\left|\frac{\zeta(2\rho)}{\zeta'(\rho)}\right|^{2k}
  \gg_k T(\log T)^{(k-1)^2},
\end{equation}
for $k\ge1$. The current version of the list prints the exponent in \eqref{eq:Ng-higher-target} as $(k+1)^2$ rather than $(k-1)^2$, which is a typographical error. We prove that it holds for every rational $k\ge1$.

\begin{theorem}\label{thm:higher-rational}
Assume the Riemann Hypothesis and that the non-trivial zeros of $\zeta(s)$ are simple. For every fixed rational number $k\ge1$,
\[
  \sum_{0<\gamma\le T}
  \left|\frac{\zeta(2\rho)}{\zeta'(\rho)}\right|^{2k}
  \gg_k T(\log T)^{(k-1)^2}.
\]
\end{theorem}

For arbitrary real exponents the same argument, followed by interpolation from below, loses an arbitrarily small power of the logarithm. We expect that this loss can be removed by adapting the method of Gao and Zhao \cite{GaoZhao}, but we do not pursue this here (see Remark~\ref{rem:rationality}).

\begin{corollary}\label{cor:higher-real-epsilon}
Assume the Riemann Hypothesis and that the non-trivial zeros of $\zeta(s)$ are simple. For every fixed real $k\ge1$ and every $\varepsilon>0$,
\[
  \sum_{0<\gamma\le T}
  \left|\frac{\zeta(2\rho)}{\zeta'(\rho)}\right|^{2k}
  \gg_{k,\varepsilon} T(\log T)^{(k-1)^2-\varepsilon}.
\]
\end{corollary}

When $1\le k<3/2$, the logarithmic exponent in Theorem~\ref{thm:higher-rational} agrees with Ng's asymptotic prediction. For $k\ge3/2$, formula \eqref{eq:Ng-higher-target} remains the corrected lower-bound target intended in the CPNT problem, but we do not claim that it gives the true order of the moment. Negative moments of $\zeta'(\rho)$ are expected to undergo a transition at $k=3/2$. The random-matrix formulae of \cite{HKO} have a phase change there. Gonek has suggested, moreover, that $|\zeta'(\rho)|^{-1}\gg\gamma^{1/3-\varepsilon}$ for infinitely many zeros, which would already defeat the conjectured order for $k>3/2$ \cite{HeapLiZhao,GaoZhao}. On the upper-bound side, item number~4 of the 2024 CPNT Problem List \cite{CPNT_ProblemList} reports the absence of conditional upper bounds for the full negative-moment sums of $\zeta'(\rho)$ discussed there. The first non-trivial upper bounds of Bui, Florea and Milinovich \cite{BuiFloreaMilinovich} concern sums restricted to a subset of the zeros.

\begin{remark}\label{rem:priority}
Theorem~\ref{thm:NgCon} was previously established in the PhD thesis of Saloni Sinha \cite{sinha2025}, and we acknowledge Sinha's priority. The proof in this paper was obtained independently. The two arguments have a similar initial contour set-up, but they differ at the left vertical edge. Sinha shifts the contour to $\Re s=1/4$, whereas we retain the line $\Re s=1-c$ and use the oscillation of $\chi(2s)/\chi(s)$. Theorem~\ref{thm:higher-rational} and Corollary~\ref{cor:higher-real-epsilon}, by contrast, are new.
\end{remark}

The moments studied here are negative moments of $\zeta'(\rho)$ in disguise, and they force the simplicity of every zero for the sum to be defined at all, and the transition at $k=3/2$ mentioned above means the conjectured order is not expected to persist.

We now describe the method, which is common to all theorems. They are proved by the mollified H\"older (or Cauchy) inequality of Rudnick and Soundararajan \cite{RudnickSound}, in the discrete form used by Milinovich and Ng \cite{MilinovichNg,MilinovichNgHigher}. In that form a first moment of $\zeta(2\rho)/\zeta'(\rho)$ against a short Dirichlet polynomial is combined with a mollified second moment. The residues $\zeta(2\rho)/\zeta'(\rho)$ are generated by the ratio
\[
\mathcal R(s)=\frac{\zeta(2s)}{\zeta(s)},
\]
so the natural mollifiers are truncations of the Dirichlet series of $\mathcal R(s)$ and of its fractional powers, whose coefficients carry the completely multiplicative Liouville sign $\lambda(n)$. For the endpoint $k=1$, that is for Theorem~\ref{thm:NgCon}, the mollifier is a smoothed truncation of $\mathcal R(s)$ itself. For the higher moments with $k=a/b$ it is the polynomial with coefficients $\lambda(n)d_{1/b}(n)$, whose integer powers keep a common Liouville sign. That sign cancels on the diagonal of the mollified first moment, leaving a non-negative sum of size $(\log T)^{k^2}$ by Lemma~\ref{lem:dk-square}.

At the endpoint the mollifier is the smoothed truncation $\mathcal L^\dagger_\xi(s)$ of $\mathcal R(s)$ of length $\xi=T^\vartheta$. Cauchy's inequality then bounds the moment below by $|M_1|^2/M_2$, where $M_1$ is the mollified first moment and $M_2$ the mollified second moment. Both are evaluated by the residue theorem on a rectangle of height $T$ with vertical sides at $\Re s=c$ and $\Re s=1-c$. The right edge carries the main term in each case. Optimising the smoothing polynomial then produces the constant $1/(4\pi)$ in Theorem~\ref{thm:NgCon}, exactly half of the conjectured $1/(2\pi)$, the factor of two being the usual loss in the Cauchy step.

For the higher moments the mollifier is the unsmoothed polynomial $A(s)$ with coefficients $\lambda(n)d_{1/b}(n)$, and Cauchy's inequality is replaced by H\"older's; the details are set out at the start of Section~\ref{sect:higher}.

The new difficulty, compared with the reciprocal moments treated in \cite{MilinovichNg,MilinovichNgHigher}, is common to all $k\ge1$ and appears on the left vertical contour. There the numerator $\zeta(2s)$ cancels the decay that $1/\zeta(s)$ alone would provide. On the line $\Re s=1-c$ the integrand has constant size up to powers of $\log T$, and the trivial bound exceeds the main term. Reflecting both zeta factors by the functional equation converts the problem into bounding integrals $\int \mathcal G(t)R^{it}\,dt$, where $\mathcal G(t)=\chi(2s)/\chi(s)$ is genuinely oscillatory. A weighted van der Corput estimate shows that each such integral is $O(\sqrt T)$ uniformly in the frequency $R>0$, which makes the left contour negligible.

The horizontal contours are controlled on a common sequence of good ordinates. That device and the oscillatory estimate, together with the Montgomery--Vaughan mean value theorem, are used for all $k\ge1$. The Selberg--Delange estimates for $d_k(n)^2$ are needed only for the higher moments. All of them are collected in Section~\ref{sect:prelim}.

\section{Notation and preliminary lemmas}\label{sect:prelim}

Throughout, $T$ is large and $c=1+1/\log T$. We write $\lambda(n)=(-1)^{\Omega(n)}$ for the Liouville function, where $\Omega(n)$ is the number of prime factors of $n$ counted with multiplicity, so that $\lambda(1)=1$ and $\lambda(n)$ is completely multiplicative. We write $\mu(n)$ for the M\"obius function, $\Lambda(n)$ for the von Mangoldt function, and $\tau(n)$ for the ordinary divisor function. We write $N(t)$ for the number of non-trivial zeros with $0<\gamma\le t$, so that the Riemann--von Mangoldt formula gives $N(t)=\frac{t}{2\pi}\log\frac t{2\pi}-\frac t{2\pi}+O(\log t)$. For $\Re s>1$,
\begin{equation}\label{eq:liouville_dirichlet}
\sum_{n=1}^\infty \frac{\lambda(n)}{n^s} = \frac{\zeta(2s)}{\zeta(s)} = \mathcal R(s).
\end{equation}
At a simple non-trivial zero $\rho$ of $\zeta(s)$, the function $\mathcal R(s)$ has a simple pole with residue $\zeta(2\rho)/\zeta'(\rho)$, whose squared modulus is the summand in Conjecture~\ref{conj:ng}. We write $\chi(s)$ for the factor in the functional equation $\zeta(s)=\chi(s)\zeta(1-s)$.

\subsection{Good ordinates}
We use one common sequence of good ordinates throughout the paper; the construction is the standard one underlying the contour arguments of Titchmarsh \cite[Theorem~14.16]{Titchmarsh} and Davenport \cite[Chapter~17]{Davenport}.

\begin{lemma}\label{lem:common-height}
Assume the Riemann Hypothesis. There are absolute constants $A_0>0$ and $\eta>0$, and a sequence $\mathcal T=\{\tau_m\}_{m\ge3}$, the index starting at $3$ so that $\log\log\tau_m$ is defined and positive, with
\[
 m<\tau_m\le m+1,
\]
such that
\begin{align}
 |\zeta(\sigma+i\tau_m)|^{-1}
 &\ll \exp\!\left(\frac{A_0\log \tau_m}{\log\log \tau_m}\right)
 &&\left(-1\le \sigma\le2\right),\label{eq:good-recip}\\
 \min_{\gamma}|\tau_m-\gamma|
 &\ge \frac{\eta}{\log\tau_m},\label{eq:good-sep}\\
 \frac{\zeta'(\sigma+i\tau_m)}{\zeta(\sigma+i\tau_m)}
 &\ll (\log\tau_m)^2
 &&(-1\le\sigma\le2).\label{eq:good-logder}
\end{align}
The implied constants are absolute.
\end{lemma}

\begin{proof}
We make the simultaneous choice explicit. For sufficiently large $m$, put $\delta=1/\log\log m$ and
\[
 F_m(t)=\sum_{|t-\gamma|\le\delta}\log\frac{2\delta}{|t-\gamma|},\qquad m<t<m+1.
\]
The proof of \cite[Theorem~14.16]{Titchmarsh} gives, uniformly for $\tfrac12\le\sigma\le2$,
\[
 \log|\zeta(\sigma+it)|\ge-F_m(t)-C\delta\log m.
\]
Moreover, $F_m\ge0$ and
\[
 \int_m^{m+1}F_m(t)\,dt
 \le 2\delta(1+\log2)\{N(m+2)-N(m-1)\}
 \le C_1\delta\log m.
\]
Consequently the set where $F_m(t)>4C_1\delta\log m$ has measure at most $1/4$. Removing also the intervals of radius $\eta/\log m$ about zeros in $(m-1,m+2]$ deletes measure at most $1/4$, if the absolute constant $\eta>0$ is sufficiently small. Choose $\tau_m$ in the remaining set. It satisfies both \eqref{eq:good-recip} on $[\tfrac12,2]$ and \eqref{eq:good-sep}, after enlarging $A_0$. The finitely many smaller $m$ are covered by decreasing $\eta$ and increasing the implied constant.

We extend \eqref{eq:good-recip} to $-1\le\sigma<\tfrac12$ by reflection. There the functional equation gives
\[
 \frac1{\zeta(\sigma+it)}
 =\frac{\chi(\sigma+it)^{-1}}{\zeta(1-\sigma-it)},
\]
and $1-\sigma$ lies in $(\tfrac12,2]$, so the bound just proved applies to the denominator. Stirling's formula yields $|\chi(\sigma+it)|^{-1}\ll t^{\sigma-1/2}\ll1$ for $\sigma<\tfrac12$, which gives \eqref{eq:good-recip} on the whole of $-1\le\sigma\le2$. Finally, the standard partial-fraction formula (see \cite[Chapter~15]{Davenport}), uniformly for $-1\le\sigma\le2$, is
\[
 \frac{\zeta'(\sigma+it)}{\zeta(\sigma+it)}
 =\sum_{|\gamma-t|\le1}\frac1{\sigma+it-\rho}+O(\log t).
\]
There are $O(\log t)$ terms in the sum, and \eqref{eq:good-sep} bounds every denominator below, since
\[
 |\sigma+it-\rho|\ge|t-\gamma|\ge\frac{\eta}{\log t}.
\]
This lower bound is independent of $\sigma$. This proves \eqref{eq:good-logder} on the same sequence, and on the whole strip $-1\le\sigma\le2$ rather than on a strip shrinking with $\tau_m$, which is what lets us apply it on the full vertical range of the contours below.
\end{proof}

\subsection{Mean values of Dirichlet polynomials}
The right vertical contours are evaluated with the mean value theorem of Montgomery and Vaughan \cite{MV}, in the convenient bilinear form of Tsang \cite[Lemma~1]{Tsang}. The mean-square case $b_n=\overline{a_n}$ is the theorem of Montgomery and Vaughan, and Tsang showed that the bilinear form below follows from it.

\begin{lemma}[Montgomery--Vaughan]\label{lem:MV}
Let $\{a_n\}_{n=1}^{\infty}$ and $\{b_n\}_{n=1}^{\infty}$ be sequences of complex numbers. Assume that $\sum_{n=1}^{\infty}n|a_n|^2<\infty$ and $\sum_{n=1}^{\infty}n|b_n|^2<\infty$, and interpret the Dirichlet series as mean-square limits of their partial sums on bounded intervals. Then for any real number $T > 0$,
\begin{equation*}
\begin{aligned}
&\int_0^T \left( \sum_{n=1}^\infty a_n n^{-it} \right)
\left( \sum_{n=1}^\infty b_n n^{it} \right)\,dt \\
&\qquad = T \sum_{n=1}^\infty a_n b_n
+ O \left( \left( \sum_{n=1}^\infty n |a_n|^2 \right)^{1/2}
\left( \sum_{n=1}^\infty n |b_n|^2 \right)^{1/2} \right).
\end{aligned}
\end{equation*}
\end{lemma}

\subsection{The oscillatory left edge}
The next two lemmas are the key to the left vertical contours. The first is a weighted form of the second-derivative test.

\begin{lemma}\label{lem:weighted-vdc}
Let $F\in C^2([U_1,U_2])$ be real-valued, suppose that $F''$ has constant sign, and suppose that $|F''(t)|\ge\kappa>0$ throughout $[U_1,U_2]$. If $g$ has bounded variation on $[U_1,U_2]$, then
\[
 \left|\int_{U_1}^{U_2} g(t)e^{iF(t)}\,dt\right|
 \ll \kappa^{-1/2}\left(\sup_{U_1\le t\le U_2}|g(t)|+\operatorname{Var}_{[U_1,U_2]}g\right).
\]
The implied constant is absolute.
\end{lemma}

\begin{proof}
The unweighted second-derivative estimate gives
\[
 \sup_{U_1\le u<v\le U_2}\left|\int_u^v e^{iF(t)}\,dt\right|
 \ll\kappa^{-1/2};
\]
see \cite[Lemma~4.4]{Titchmarsh}. Partial summation in the Riemann--Stieltjes sense then gives the displayed weighted form. Titchmarsh's own weighted version, \cite[Lemma~4.5]{Titchmarsh}, assumes instead that $g/F'$ is monotonic, which the weight $|\mathcal G|$ of Lemma~\ref{lem:osc-edge} need not satisfy; that is why the bounded-variation form is derived here.
\end{proof}

Reflecting $\mathcal R(s)$ on the line $\Re s=1-c$ produces the quotient
\begin{equation}\label{eq:Gdef}
\mathcal G(t)=\frac{\chi(2(1-c)+2it)}{\chi((1-c)+it)},
\end{equation}
and every left-edge integral in this paper reduces to the following estimate.

\begin{lemma}\label{lem:osc-edge}
There is an absolute constant $t_0\ge1$ with the following property. Let $T$ be large and let $\mathcal G(t)$ be as in \eqref{eq:Gdef}. Then, uniformly for every real $R>0$ and every $t_0\le U_1<U_2\le2T$,
\[
 \int_{U_1}^{U_2}\mathcal G(t)\,R^{it}\,dt\ll\sqrt T.
\]
\end{lemma}

\begin{proof}
For $t\ge t_0$ the function $\chi(2s)/\chi(s)$ with $s=(1-c)+it$ is analytic and non-vanishing, the zeros of $\sin(\pi s/2)$ and the poles of $\Gamma(1-s)$ being real. Hence $\mathcal G(t)=|\mathcal G(t)|e^{i\arg \mathcal G(t)}$ with a smooth choice of argument. By Stirling's formula, for $t\ge t_0$,
\[
|\mathcal G(t)| \asymp 2^{2c-3/2}\left(\frac{t}{2\pi}\right)^{c-1}
\qquad\text{and}\qquad
\frac{d}{dt}|\mathcal G(t)|\ll|\mathcal G(t)|\left(\frac1{t\log T}+\frac1{t^2}\right).
\]
Since $t^{c-1}=t^{1/\log T}\ll1$ for $t\le2T$, the modulus $|\mathcal G|$ is bounded on $[t_0,2T]$. Integrating the derivative bound, and using
\[
 \int_{t_0}^{2T}\left(\frac1{t\log T}+\frac1{t^2}\right)dt\ll1,
\]
gives $\operatorname{Var}_{[t_0,2T]}|\mathcal G|\ll1$.

We turn to the phase. Stirling's asymptotic series for $\log\Gamma$ gives the expansion $\arg\chi(\sigma+iu) = -u\log(u/2\pi)+u+\tfrac\pi4+O(1/u)$, uniformly for bounded $\sigma$; see \cite[\S4.12]{Titchmarsh}. We use Stirling's formula and its differentiated versions to control the phase derivatives. Put $\psi_R(t)=\arg \mathcal G(t)+t\log R$. These estimates give
\[
  \psi_R'(t) = -\log\frac{2t}{\pi R} + O\!\Big(\frac1t\Big), \qquad \psi_R''(t) = -\frac1t + O\!\Big(\frac1{t^2}\Big).
\]
Three consequences follow, once $t_0$ is a sufficiently large absolute constant. First, $\psi_R''$ has constant sign on $[t_0,2T]$. Second, that sign is independent of $R$. Third, $|\psi_R''(t)|\gg T^{-1}$ throughout $[t_0,2T]$. Writing $\mathcal G(t)R^{it}=|\mathcal G(t)|e^{i\psi_R(t)}$, Lemma~\ref{lem:weighted-vdc} applies with $g=|\mathcal G|$ and $F=\psi_R$, and yields the claim.
\end{proof}

\subsection{Divisor sums}
Finally, for the higher moments we need standard consequences of Selberg--Delange theory. For $z>0$, let $d_z(n)$ be defined by
\[
  \zeta(s)^z=\sum_{n=1}^{\infty}\frac{d_z(n)}{n^s},\qquad \Re s>1.
\]
Since $\mathcal R(s)=\prod_p(1+p^{-s})^{-1}$, one has
\begin{equation}\label{eq:Rz-coeff}
  \mathcal R(s)^z=\sum_{n=1}^{\infty}\frac{\lambda(n)d_z(n)}{n^s},
  \qquad \Re s>1.
\end{equation}
\begin{lemma}\label{lem:dk-square}
Let $k\ge1$ be fixed. There is a function $F_k(s)$, holomorphic for $\Re s>\tfrac12$ and non-vanishing in a neighbourhood of $s=1$, such that
\begin{equation}\label{eq:dk-factorisation}
 \sum_{n=1}^{\infty}\frac{d_k(n)^2}{n^s}=\zeta(s)^{k^2}F_k(s)
 \qquad(\Re s>1).
\end{equation}
Moreover, uniformly for $y\ge3$,
\[
 \sum_{n\le y}\frac{d_k(n)^2}{n}\asymp_k(\log y)^{k^2},
\]
and, for $j=0,1,2$ and $0<\delta\le1$,
\[
 \sum_{n=1}^{\infty}\frac{d_k(n)^2(\log n)^j}{n^{1+\delta}}
 \ll_k \delta^{-k^2-j}.
\]
\end{lemma}

\begin{proof}
At a prime $p$ the local factor of the series on the left of \eqref{eq:dk-factorisation} is
\[
 \sum_{r=0}^{\infty}d_k(p^r)^2p^{-rs}
 =1+k^2p^{-s}+O_k(p^{-2\Re s}),
\]
whereas the local factor of $\zeta(s)^{k^2}$ is
$1+k^2p^{-s}+O_k(p^{-2\Re s})$. Their quotient therefore defines an Euler product that converges absolutely and uniformly in $\Re s\ge\tfrac12+\delta_0$, for each fixed $\delta_0>0$. This gives a holomorphic and bounded $F_k$ in each such half-plane. Each local factor is non-zero for real $s\in(\tfrac12,2]$ and hence in a complex neighbourhood, while the uniformly convergent tail of the product is non-vanishing, so $F_k$ is non-zero in a neighbourhood of $s=1$. In particular, on the real interval $1<s\le2$, the quantities $F_k(s)$ and $F_k(s)^{-1}$ are bounded in terms of $k$.

The Selberg--Delange theorem requires the correcting factor to be holomorphic and bounded in a region $\sigma>1-c_0/\log(|t|+2)$. Taking $c_0$ small, any such region is contained in a half-plane $\Re s\ge\tfrac12+\delta_0$, so the hypotheses are met. The theorem applied to \eqref{eq:dk-factorisation} gives
\[
 \sum_{n\le y}d_k(n)^2
 =\kappa_k\, y(\log y)^{k^2-1}+O_k\bigl(y(\log y)^{k^2-2}\bigr)
\]
with $\kappa_k>0$. Partial summation yields the first displayed estimate. Finally, evaluate \eqref{eq:dk-factorisation} and its first two derivatives at the real point $s=1+\delta$. Using $\zeta^{(r)}(1+\delta)\ll_r\delta^{-r-1}$ and the boundedness of the derivatives of $F_k$ on $1<s\le2$, this gives the remaining estimates. See also \cite[Chapter~II.5]{Tenenbaum}.
\end{proof}

\section{The explicit second moment: set-up}\label{sect:setup-second}
We follow the strategy of Milinovich--Ng \cite{MilinovichNg}, which is based on the work of Rudnick and Soundararajan \cite{RudnickSound}.

Let $\xi = T^\vartheta$ where $0 < \vartheta < 1$ is fixed. Let $P(x)$ be a polynomial with real coefficients satisfying $P(0)=0$ and $P(1)=1$, and define the Dirichlet polynomial
\begin{equation*}
\mathcal{L}_\xi^\dagger(s) = \sum_{n \le \xi} \frac{\lambda(n)}{n^s} P\left(\frac{\log(\xi/n)}{\log \xi}\right).
\end{equation*}
Throughout Sections~\ref{sect:setup-second} to~\ref{sect:explicit-proof} we abbreviate the smoothing weight by
\begin{equation}\label{eq:Pn-def}
 P_n = P\left(\frac{\log(\xi/n)}{\log \xi}\right)
 \qquad(n\le\xi).
\end{equation}
The dagger records that the truncation is smoothed. By \eqref{eq:liouville_dirichlet}, $\mathcal{L}_\xi^\dagger(s)$ is a smoothed Dirichlet polynomial approximating $\mathcal R(s)$. The squared moduli of the residues of $\mathcal R(s)$ at the non-trivial zeros are precisely the summands appearing in Conjecture~\ref{conj:ng}. This motivates its use as a mollifier and explains the appearance of $\mathcal R(s)$ in the analysis of $M_1$ below.

Assuming the Riemann Hypothesis, for any non-trivial zero $\rho = 1/2 + i\gamma$ of $\zeta(s)$ we have $\overline{\mathcal{L}_\xi^\dagger(\rho)} = \mathcal{L}_\xi^\dagger(1-\rho)$; here the reality of the coefficients of $P$ is used. From this observation and Cauchy's inequality, it follows that
\begin{equation} \label{eq:ratio}
\sum_{0 < \gamma \le T} \left| \frac{\zeta(2\rho)}{\zeta'(\rho)} \right|^2 \ge \frac{|M_1|^2}{M_2}
\end{equation}
where
\[
M_1 = \sum_{0 < \gamma \le T} \frac{\zeta(2\rho)}{\zeta'(\rho)} \mathcal{L}_\xi^\dagger(1-\rho) \quad \text{and} \quad M_2 = \sum_{0 < \gamma \le T} |\mathcal{L}_\xi^\dagger(\rho)|^2.
\]
The inequality is read trivially if $M_2=0$. In fact $M_2>0$ for large $T\in\mathcal T$: the bracketed factor in the main term of Proposition~\ref{prop:M2} is at least $\vartheta\int_0^1P(x)^2\,dx$, which is positive because $P$ is real and $P(1)=1$. That proposition is proved independently of \eqref{eq:ratio}.

\begin{remark}\label{rem:roles-of-hypotheses}
The two hypotheses play distinct roles. The simplicity of the zeros is used to guarantee that the sum in Conjecture~\ref{conj:ng} is well-defined, that is, that $\zeta'(\rho)\ne0$. It is used again to ensure that the poles of $\mathcal R(s)$ at the non-trivial zeros are simple, with residue $\zeta(2\rho)/\zeta'(\rho)$, in the residue computation of Section~\ref{sect:M1}. Proposition~\ref{prop:M2} needs no simplicity: the residue of $\zeta'(s)/\zeta(s)$ at a zero of multiplicity $m_\rho$ is $m_\rho$, which matches the multiplicity-counted sum $M_2$. The Riemann Hypothesis is used to place all non-trivial zeros on the critical line, which yields the symmetry $\overline{\mathcal{L}_\xi^\dagger(\rho)} = \mathcal{L}_\xi^\dagger(1-\rho)$ underlying \eqref{eq:ratio}. It is used again in the zero-free region and size estimates employed in Sections~\ref{sect:M1} and~\ref{sect:M2}.
\end{remark}

To prove Theorem \ref{thm:NgCon}, we must evaluate the asymptotics of $M_1$ and $M_2$. In Sections~\ref{sect:M1} and~\ref{sect:M2}, respectively, we establish the following estimates on the good sequence $\mathcal T$ of Lemma~\ref{lem:common-height}.

\begin{proposition} \label{prop:M1}
Assume the Riemann Hypothesis, let $0 < \vartheta < 1$ be fixed, and assume the zeros of $\zeta(s)$ are all simple. For $T$ in the sequence $\mathcal T$ of Lemma~\ref{lem:common-height},
\begin{equation*}
M_1 = \frac{T}{2\pi} (\log T) \; \vartheta \int_0^1 P(x)\,dx + O(T).
\end{equation*}
\end{proposition}

\begin{proposition} \label{prop:M2}
Assume the Riemann Hypothesis and let $0 < \vartheta < 1$ be fixed. For $T$ in the sequence $\mathcal T$ of Lemma~\ref{lem:common-height},
\begin{equation*}
M_2 = \frac{T}{2\pi}(\log T)^2 \left[ \vartheta \int_0^1 P(x)^2\,dx + \vartheta^2 \left( \int_0^1 P(x)\,dx \right)^2 \right] + O(T \log T).
\end{equation*}
\end{proposition}

In Section~\ref{sect:explicit-proof} we substitute these asymptotics into \eqref{eq:ratio} and optimise the polynomial $P(x)$ to prove Theorem~\ref{thm:NgCon}.

\section{Proof of Proposition \ref{prop:M1}}\label{sect:M1}
Let $T\in\mathcal T$ with $T\ge4$, so that $\log\log T>0$, and set $c = 1 + 1/\log T$. By Lemma~\ref{lem:common-height} there is a fixed constant $A_0>0$ such that
\begin{equation} \label{eq:zeta_bound}
|\zeta(\sigma+iT)|^{-1} \ll \exp\left(\frac{A_0 \log T}{\log \log T}\right)
\end{equation}
uniformly for $-1 \le \sigma \le 2$. Assume that all the non-trivial zeros of $\zeta(s)$ are simple. Then the poles of $\mathcal R(s)$ inside the rectangle below are exactly the zeros $\rho$ with $1<\gamma\le T$, each simple and with residue $\zeta(2\rho)/\zeta'(\rho)$. Indeed, the pole of $\zeta(2s)$ at $s=\tfrac12$ lies on the real axis, below the bottom edge, and there are no zeros with $0<\gamma\le1$, the least positive ordinate exceeding $14$. Hence the residue theorem gives
\begin{align*}
    M_1 &= \frac{1}{2\pi i} \left( \int_{c+i}^{c+iT} + \int_{c+iT}^{1-c+iT} + \int_{1-c+iT}^{1-c+i} + \int_{1-c+i}^{c+i} \right) \mathcal R(s)\, \mathcal{L}_\xi^\dagger(1-s) \,ds \\
    &= I_1 + I_2 + I_3 + I_4,
\end{align*}
say. We will show that the main contribution to $M_1$ comes from the right vertical contour $I_1$, and that the integrals over the other contours only contribute to the error term.

\subsection*{The right edge}
By \eqref{eq:liouville_dirichlet}, the Dirichlet series for $\mathcal R(s)$ is generated by the Liouville function. In the notation \eqref{eq:Pn-def} we may rewrite $I_1$ as
\begin{equation*}
I_1 = \frac{1}{2\pi} \int_1^T \sum_{m=1}^\infty \frac{\lambda(m)}{m^{c+it}} \sum_{n \le \xi} \frac{\lambda(n) P_n}{n^{1-c-it}} \,dt.
\end{equation*}
Apply Lemma \ref{lem:MV} at the heights $T$ and $1$ and difference, with $a_m = \lambda(m)m^{-c}$ and $b_n = \lambda(n) P_n n^{-1+c}$. This gives
\begin{equation*}
I_1 = \frac{(T-1)}{2\pi} \sum_{n \le \xi} \frac{\lambda(n)^2 P_n}{n} + O\left( \left( \sum_{m=1}^\infty \frac{\lambda(m)^2}{m^{2c-1}} \right)^{1/2} \left( \sum_{n \le \xi} \frac{\lambda(n)^2 P_n^2}{n^{1-2c}} \right)^{1/2} \right),
\end{equation*}
where the error term is
\begin{equation*}
O \left(\zeta(2c-1)^{1/2} \left( \xi^{2c} \right)^{1/2} \right) \ll \xi \sqrt{\log T},
\end{equation*}
recalling that $c = 1 + 1/\log T$, so that $\zeta(2c-1) \asymp \log T$ and $\xi^c \asymp \xi$. Since $\lambda(n)^2 = 1$ for all integers $n$, partial summation yields
\begin{equation*}
\sum_{n \le \xi} \frac{1}{n} P\left(\frac{\log(\xi/n)}{\log \xi}\right) = \log \xi \int_0^1 P(x) \,dx + O(1).
\end{equation*}
Hence
\begin{equation}\label{eq:I1}
I_1 = \frac{T}{2\pi} (\log \xi) \int_0^1 P(x) \,dx + O(\xi \sqrt{\log T} + T).
\end{equation}

\subsection*{The horizontal edges}
To estimate the upper horizontal contour $I_2$, we use the convexity bound for $\zeta(s)$ (see \cite[Chapter~V]{Titchmarsh}) and the uniform elementary bound
\[
 |\mathcal{L}_\xi^\dagger(1-\sigma-iT)|
 \ll (1+\log\xi)\xi^{\max(\sigma,0)}.
\]
Since $[1-c,c]\subset[-1,2]$, the reciprocal bound \eqref{eq:zeta_bound} holds throughout the horizontal segment. We therefore have
\begin{equation*}
I_2 \ll (1+\log\xi)\exp\left(\frac{A_0 \log T}{\log \log T}\right) \int_{1-c}^c |\zeta(2\sigma+2iT)| \xi^\sigma \,d\sigma.
\end{equation*}
For $\sigma \ge 1/2$ the integrand is $\ll \xi^\sigma \log T \le \xi \log T$. For $0\le\sigma<1/2$ the convexity bound, in the strip $0\le\Re(2s)\le1$, gives $|\zeta(2\sigma+2iT)| \ll T^{1/2 - \sigma + o(1)}$. For $1-c\le\sigma<0$ reflection gives $|\zeta(2\sigma+2iT)|\ll T^{1/2-2\sigma}\ll T^{1/2-\sigma+o(1)}$, since $-\sigma\le1/\log T$ makes $T^{-\sigma}$ bounded. In both of these ranges the integrand is therefore $\ll T^{1/2 - \sigma(1-\vartheta) + o(1)} \le T^{1/2 + o(1)}$. Hence
\begin{equation}\label{eq:I2}
I_2 \ll (1+\log\xi)\left( \xi + T^{1/2} \right) \exp\left(\frac{A_0 \log T}{\log \log T}\right) = o(T).
\end{equation}
The bottom path $I_4$ has bounded length. As $T$ varies, these paths lie in a fixed compact subset of the line $\Im s=1$ on which neither $\zeta(s)$ nor $\zeta(2s)$ has a zero or pole, so $\mathcal R(s)$ is uniformly bounded there. Since $|\mathcal L_\xi^\dagger(1-s)|\ll(1+\log\xi)\xi^{\max(\Re s,0)}$, it follows that $I_4\ll\xi(1+\log\xi)=o(T)$.

\subsection*{The left edge}
The left vertical contour $I_3$ cannot be bounded trivially. On $\Re s = 1-c$ one has $2\Re s = -2/\log T < 0$, so the functional equation gives $|\zeta(2s)| \ll |t|^{1/2}\log T$, while $|\zeta(s)|^{-1} \ll |t|^{-1/2}\log T$. Their product is $\ll(\log T)^2$ and does not decay in $t$. Since $\mathcal L^\dagger_\xi(1-s)\ll\log\xi$ on that line, the trivial bound would give $I_3 \ll T(\log T)^3$, which exceeds the main term $\asymp T\log T$. In contrast to Milinovich--Ng \cite{MilinovichNg}, where the numerator is merely $1/\zeta(s)$ and decays like $|t|^{-1/2}$, here the extra factor $\zeta(2s)$ cancels that decay. We must instead exploit the oscillation of $\chi(2s)/\chi(s)$.

\begin{remark}\label{rem:sinha-divergence}
It is precisely at this stage of the argument that the present approach diverges from the independent proof given by Sinha \cite{sinha2025}. To handle the contribution from the $\zeta(2s)$ factor, Sinha evaluates the corresponding integral by shifting the line of integration further left to $\Re s=\tfrac14$. We instead keep the contour at $\Re s=1-c$. After reflecting through the functional equation, the problem becomes an integral governed by the highly oscillatory factor $\chi(2s)/\chi(s)$, which Lemma~\ref{lem:osc-edge} controls.
\end{remark}

Applying the functional equation $\zeta(w)=\chi(w)\zeta(1-w)$ to both $w=2s$ and $w=s$ gives
\begin{equation*}
  \mathcal R(s) \;=\; \frac{\chi(2s)}{\chi(s)}\,\frac{\zeta(1-2s)}{\zeta(1-s)} \;=\; \mathcal G(t)\,\frac{\zeta(1-2s)}{\zeta(1-s)}, \qquad s = (1-c)+it,
\end{equation*}
with $\mathcal G(t)$ as in \eqref{eq:Gdef}. Since $\Re(1-2s)=2c-1>1$ and $\Re(1-s)=c>1$, both reflected zeta factors lie in the region of absolute convergence. Expanding, we have
\[
\frac{\zeta(1-2s)}{\zeta(1-s)}\, \mathcal{L}_\xi^\dagger(1-s) = \sum_{n=1}^{\infty}\beta_n n^{it}
\]
with
\[
\beta_n=\sum_{a^2 b m=n,\, m\le\xi} a^{-(2c-1)}\mu(b)b^{-c}\lambda(m)P_m m^{-c}.
\]

The factor $m^{-c}$ comes from $\mathcal L_\xi^\dagger(1-s)$ and is required in the $\ell^1$ estimate below. The series converges absolutely, again because $\Re(1-2s)=2c-1>1$ and $\Re(1-s)=c>1$, so it may be integrated term by term against the bounded function $\mathcal G(t)$. Thus
\begin{equation*}
  I_3 \;=\; -\frac{1}{2\pi}\sum_{n=1}^{\infty}\beta_n\, \mathcal J_n,
\end{equation*}
where $\mathcal J_n = \int_1^T \mathcal G(t)\,n^{it}\,dt$. The initial interval $[1,t_0]$ contributes $O(1)$ to each $\mathcal J_n$, since $\mathcal G$ is continuous and non-vanishing on the fixed compact segment $\{(1-c)+it:1\le t\le t_0\}$, uniformly in $T$. Lemma~\ref{lem:osc-edge} with $R=n$ bounds the rest, so that, uniformly in $n$,
\begin{equation*}
  |\mathcal J_n| \;\ll\; \sqrt{T}.
\end{equation*}
Moreover,
\[
\sum_{n}|\beta_n| \le \zeta(2c-1)\,\frac{\zeta(c)}{\zeta(2c)}\,\Big(\max_{[0,1]}|P|\Big)\zeta(c) \ll (\log T)^3,
\]
since $\zeta(c)\ll\log T$ and $\zeta(2c-1)\ll\log T$. Therefore
\begin{equation}\label{eq:I3}
  I_3 \;\ll\; \Big(\sum_{n}|\beta_n|\Big)\,\sup_n|\mathcal J_n| \;\ll\; \sqrt{T}\,(\log T)^3 \;=\; o(T).
\end{equation}

\begin{remark}\label{rem:multiplicative-fe}
The companion left edge $J_3$, in the proof of Proposition~\ref{prop:M2}, is governed by the additive functional equation for $\zeta'(s)/\zeta(s)$. Its non-oscillatory remainder $-\log(t/2\pi)$ is what produces the extra main term $K$ there, and $K_B$ in Section~\ref{sect:higher}. The functional equation for $\mathcal R(s)$ acts multiplicatively instead. Every term $\beta_n n^{it}$ therefore carries the factor $\mathcal G(t)$, whose phase is stationary only at the isolated point near $t=\pi n/2$. No non-oscillatory summand survives, and $I_3$ contributes no secondary main term.
\end{remark}

Combining \eqref{eq:I1}, \eqref{eq:I2}, and \eqref{eq:I3} with the bound for $I_4$ and substituting $\log \xi = \vartheta \log T$ completes the proof of Proposition \ref{prop:M1}.

\section{Proof of Proposition \ref{prop:M2}}\label{sect:M2}
We again let $T\in\mathcal T$ with $T \geq 4$ and $c = 1 +1/\log T$. Assuming the Riemann Hypothesis, we have
\[
M_2 = \sum_{0 < \gamma \leq T} \mathcal{L}_\xi^\dagger(\rho) \mathcal{L}_\xi^\dagger(1-\rho).
\]
The only poles of the integrand below that lie inside the rectangle are the non-trivial zeros $\rho$ with $1<\gamma\le T$, the residue at a zero of multiplicity $m_\rho$ being $m_\rho\mathcal L^\dagger_\xi(\rho)\mathcal L^\dagger_\xi(1-\rho)$. The pole of $\zeta'(s)/\zeta(s)$ at $s=1$ lies on the real axis, below the bottom edge, and there are no zeros with $0<\gamma\le1$. Hence the residue theorem gives
\begin{align*}
M_2 &= \frac{1}{2\pi i} \left( \int_{c+i}^{c+iT} + \int_{c+iT}^{1-c+iT} + \int_{1-c+iT}^{1-c+i} + \int_{1-c+i}^{c+i} \right) \frac{\zeta'(s)}{\zeta(s)} \mathcal{L}_\xi^\dagger(s) \mathcal{L}_\xi^\dagger(1-s) \,ds \\
&= J_1 + J_2 + J_3 + J_4,
\end{align*}
say.

Since $T\in\mathcal T$, estimate \eqref{eq:good-logder} gives
\[
\frac{\zeta'(\sigma+iT)}{\zeta(\sigma+iT)}\ll(\log T)^2
\qquad(1-c\le\sigma\le c),
\]
on the same fixed sequence as in Proposition~\ref{prop:M1}. We further use the uniform elementary bound
\[
|\mathcal{L}_\xi^\dagger(\sigma + it)| \ll (1+\log \xi)\,\xi^{\max(1-\sigma,0)},
\]
where the factor $1+\log\xi$ is needed to make the estimate valid also for $\sigma>1$. Consequently, on the horizontal contours the product of the two Dirichlet polynomials is $\ll \xi(1+\log \xi)^2$, and hence $J_2\ll \xi(\log T)^4$. On the bottom path, of bounded length at height $1$, the factor $\zeta'(s)/\zeta(s)$ is uniformly bounded, its pole lying at distance at least $1$ from the fixed compact segment, so also $J_4\ll\xi(\log T)^4$.

We relate the left vertical contour $J_3$ to $J_1$ via the logarithmic derivative of the functional equation, which after an application of Stirling's formula gives
\[
\frac{\zeta'(1-s)}{\zeta(1-s)} = -\frac{\zeta'(s)}{\zeta(s)} - \log\frac{|t|}{2\pi} + O\left(\frac{1}{|t|}\right)
\]
uniformly for $1 \leq |t| \leq T$. The remainder contributes
\[
 \ll\int_1^T t^{-1}\left|\mathcal L^\dagger_\xi(c+it)\,\mathcal L^\dagger_\xi(1-c-it)\right|\,dt
 \ll\xi(\log T)^3,
\]
since $\mathcal L^\dagger_\xi(c+it)\ll\log\xi$ and $\mathcal L^\dagger_\xi(1-c-it)\ll\xi\log\xi$, and the integral of $t^{-1}$ contributes a further $\log T$. Changing variables $s \to 1-s$ and taking the complex conjugate, we find
\[
M_2 = \overline{K} + 2\Re J_1 + O(\xi(\log T)^4),
\]
where
\begin{equation*}
K = \frac{1}{2\pi} \int_1^T \log\left(\frac{t}{2\pi}\right) \mathcal{L}_\xi^\dagger(c+it) \mathcal{L}_\xi^\dagger(1-c-it) \,dt.
\end{equation*}
The main term in $K$ is real, so $\overline K$ has the same asymptotic as $K$.

To evaluate $K$, set
\[
 H(t) = \int_1^t \mathcal{L}_\xi^\dagger(c+iu) \mathcal{L}_\xi^\dagger(1-c-iu) \,du.
\]
The two factors have coefficient exponents summing to $1$. Lemma \ref{lem:MV}, applied at the heights $t$ and $1$ and differencing, therefore gives, uniformly in $t$,
\begin{equation} \label{eq:At}
H(t) = t \sum_{n \le \xi} \frac{\lambda(n)^2 P_n^2}{n} + O(\xi \sqrt{\log T}).
\end{equation}
Partial summation on the diagonal, as at \eqref{eq:I1}, converts this into
\begin{equation}\label{eq:At-2}
 H(t) = t \log \xi \int_0^1 P(x)^2 \,dx + O\big(t+\xi\sqrt{\log T}\big),
\end{equation}
the new $O(t)$ arising from the $O(1)$ error in the partial summation. Integrating by parts, we write
\begin{equation*}
K = \frac{1}{2\pi} \int_1^T \log\left(\frac{t}{2\pi}\right) H'(t) \,dt = \frac{1}{2\pi} \left[ \log\left(\frac{t}{2\pi}\right) H(t) \right]_1^T - \frac{1}{2\pi} \int_1^T \frac{H(t)}{t} \,dt.
\end{equation*}
Substituting \eqref{eq:At-2} into this expression gives
\begin{align*}
K &= \frac{T}{2\pi} \log\left(\frac{T}{2\pi}\right) \log \xi \int_0^1 P(x)^2 \,dx - \frac{T}{2\pi} \log \xi \int_0^1 P(x)^2 \,dx + O(T \log T)\\
&= \frac{T}{2\pi} \log\left(\frac{T}{2\pi}\right) \log \xi \int_0^1 P(x)^2 \,dx + O(T \log T)\\
&= \frac{\vartheta}{2\pi} T (\log T)^2 \int_0^1 P(x)^2 \,dx + O(T \log T),
\end{align*}
using $\log(T/2\pi) \log \xi = \vartheta (\log T)^2 + O(\log T)$.

For $J_1$, we consider the Dirichlet series
\begin{equation*}
-\frac{\zeta'(s)}{\zeta(s)} \mathcal{L}_\xi^\dagger(s) = \left( \sum_{h=1}^\infty \Lambda(h) h^{-s} \right) \left( \sum_{\ell \le \xi} \lambda(\ell) P_\ell\, \ell^{-s} \right) = \sum_{n=1}^\infty \alpha_n n^{-s},
\end{equation*}
where the convolution coefficients are
\[
 \alpha_n = \sum_{\substack{h\ell=n\\ \ell\le\xi}} \Lambda(h) \lambda(\ell) P_\ell .
\]
We first record a sufficient mean-square bound for these coefficients. Put $M_P=\max_{0\le x\le1}|P(x)|$. By Cauchy's inequality over the divisors of $n$,
\[
 |\alpha_n|^2
 \le M_P^2 \tau(n)\sum_{r\mid n}\Lambda(r)^2
 \le M_P^2 \tau(n)(\log n)^2;
\]
indeed, if $n=\prod p^{e_p}$, then
$\sum_{r\mid n}\Lambda(r)^2=\sum_p e_p(\log p)^2\le(\sum_p e_p\log p)^2=(\log n)^2$.
Consequently, with $\delta=c-1=1/\log T$,
\begin{align*}
 \sum_{n=1}^\infty |\alpha_n|^2n^{1-2c}
 &\ll_P\sum_{n=1}^\infty\frac{\tau(n)(\log n)^2}{n^{1+2\delta}}
 =\left.\frac{d^2}{ds^2}\zeta(s)^2\right|_{s=1+2\delta}
 \ll_P\delta^{-4}\ll_P(\log T)^4.
\end{align*}
Also $\sum_{n\le\xi}nP_n^2n^{2c-2}\ll_P\xi^2$. Hence Lemma~\ref{lem:MV}, applied at the heights $T$ and $1$ and differencing, with $a_n = \alpha_n n^{-c}$ and $b_n = \lambda(n) P_n n^{c-1}$, gives
\begin{equation*}
J_1 = -\frac{T}{2\pi} \sum_{n \le \xi} \frac{\alpha_n \lambda(n) P_n}{n} + O_P(\xi (\log T)^2),
\end{equation*}
the replacement of $T-1$ by $T$ costing $O_P((\log T)^2)$, which the displayed error absorbs.

To evaluate the main sum, substitute $\alpha_n$ and use the complete multiplicativity of the Liouville function. This gives
\begin{equation*}
\sum_{n \le \xi} \frac{\alpha_n \lambda(n) P_n}{n} = \sum_{h\ell \le \xi} \frac{\Lambda(h) \lambda(\ell) P_\ell \, \lambda(h\ell) P_{h\ell}}{h\ell} = \sum_{h\ell \le \xi} \frac{\Lambda(h) \lambda(h) \, P_\ell P_{h\ell}}{h\ell},
\end{equation*}
since $\lambda(\ell)\lambda(h\ell) = \lambda(\ell)^2 \lambda(h) = \lambda(h)$.

We carry out the sum over $h$ first, retaining the dependence of $P_{h\ell}$ on $h$. The factor $\Lambda(h)\lambda(h)$ is supported on prime powers, with $\Lambda(p)\lambda(p) = -\log p$ on primes and a convergent contribution from higher prime powers. Write $y = \frac{\log(\xi/\ell)}{\log\xi}$, so that $P_\ell = P(y)$. Mertens' theorem in the form $\sum_{p \le z} \frac{\log p}{p} = \log z + O(1)$, together with partial summation, then gives
\begin{equation*}
\sum_{h \le \xi/\ell} \frac{\Lambda(h) \lambda(h)}{h} P_{h\ell} = -\log\xi \int_0^y P(u) \,du + O(1)
\end{equation*}
uniformly in $\ell\le\xi$. The remainder is $\ll_P1$ because the partial-summation error involves only $\max_{[0,1]}|P|$ and $\operatorname{Var}_{[0,1]}P$, and the case $\xi/\ell<2$ is trivial. Summing over $\ell$ and applying partial summation once more, as in Section \ref{sect:M1}, we obtain
\begin{equation*}
\sum_{n \le \xi} \frac{\alpha_n \lambda(n) P_n}{n} = -(\log\xi)^2 \int_0^1 P(y) \left( \int_0^y P(u) \,du \right)\,dy + O(\log\xi).
\end{equation*}

Write $Q(y) = \int_0^y P(u)\,du$. The integrand equals $Q'(y)Q(y) = \tfrac{1}{2}\big(Q(y)^2\big)'$, so the double integral equals $\tfrac{1}{2} Q(1)^2 = \tfrac{1}{2}\big(\int_0^1 P\big)^2$. Hence
\begin{equation*}
\sum_{n \le \xi} \frac{\alpha_n \lambda(n) P_n}{n} = -\frac{1}{2}(\log\xi)^2 \left( \int_0^1 P(x) \,dx \right)^2 + O(\log\xi).
\end{equation*}
Substituting this into our expression for $J_1$, and recalling $\log\xi = \vartheta \log T$, gives
\begin{equation*}
2\Re J_1 = \frac{\vartheta^2}{2\pi} T (\log T)^2 \left( \int_0^1 P(x) \,dx \right)^2 + O(T \log T).
\end{equation*}

Combining the estimates for $K$ and $2\Re J_1$ establishes Proposition \ref{prop:M2}.

\section{Proof of Theorem \ref{thm:NgCon}}\label{sect:explicit-proof}
Let $T$ be large. We first record the bookkeeping around the good ordinate.

\begin{remark}\label{rem:good-height}
There is an ordinate $T'\in\mathcal T$ with $T-2<T'\le T$, because $\mathcal T$ contains a member of every interval $(j,j+1]$ with $j\ge3$. Three consequences will be used without further comment. First, the summands are non-negative, so the sum up to $T$ is at least the sum up to $T'$. Second, $T'=T+O(1)$, so replacing $T'$ by $T$ in the main terms changes them only by lower-order terms. Third, Propositions~\ref{prop:M1} and~\ref{prop:M2} are applied at the height $T'$, so that the mollifier there has length $\xi=(T')^\vartheta$; since $\log\xi=\vartheta\log T+O(1)$, this too is absorbed in the error terms.
\end{remark}

The parameters $\vartheta$ and, below, $m$ are fixed before $T$ tends to infinity. Both propositions give their main terms with relative error $O_{\vartheta,P}(1/\log T)$, provided $\int_0^1P(x)\,dx$ is non-zero, as it is for the family $R_m$ used below, and the main term of $M_2$ is bounded below by a positive multiple of $T(\log T)^2$, because $\int_0^1P^2>0$. Hence the quotient $|M_1|^2/M_2$ equals the quotient of the two main terms times $1+O_{\vartheta,P}(1/\log T)$. Combining \eqref{eq:ratio} with Propositions \ref{prop:M1} and \ref{prop:M2} therefore gives
\begin{equation*}
\sum_{0 < \gamma \le T} \left| \frac{\zeta(2\rho)}{\zeta'(\rho)} \right|^2
\ge \frac{ \big( \tfrac{\vartheta}{2\pi} T \log T \int_0^1 P \big)^2 }{ \tfrac{1}{2\pi} T (\log T)^2 \big[ \vartheta \int_0^1 P^2 + \vartheta^2 \big( \int_0^1 P \big)^2 \big] } + o(T),
\end{equation*}
which simplifies to
\begin{equation*}
\sum_{0 < \gamma \le T} \left| \frac{\zeta(2\rho)}{\zeta'(\rho)} \right|^2
\ge \frac{T}{2\pi} \cdot \frac{ \vartheta \big( \int_0^1 P \big)^2 }{ \int_0^1 P^2 + \vartheta \big( \int_0^1 P \big)^2 } + o(T).
\end{equation*}

To maximise this lower bound over all admissible polynomials $P(x)$, set $r = \big( \int_0^1 P \big)^2 \big/ \int_0^1 P^2$. Since $P$ is real, Cauchy's inequality gives
\begin{equation*}
\left( \int_0^1 P(x) \,dx \right)^2 \le \left( \int_0^1 1^2 \,dx \right) \left( \int_0^1 P(x)^2 \,dx \right) = \int_0^1 P(x)^2 \,dx,
\end{equation*}
so that $r \le 1$.

Equality would require $P$ to be constant almost everywhere, and hence cannot be attained within the admissible class $P(0)=0$, $P(1)=1$. Nevertheless, the supremum is $1$. Indeed, for every positive integer $m$ the polynomial
\[
R_m(x)=1-(1-x)^m
\]
is admissible, and
\[
\int_0^1 R_m(x)\,dx=\frac{m}{m+1},
\qquad
\int_0^1 R_m(x)^2\,dx=1-\frac{2}{m+1}+\frac{1}{2m+1},
\]
so that the corresponding ratio is
\[
 r_m=\frac{2m+1}{2m+2}\longrightarrow1
 \qquad(m\to\infty).
\]
In terms of $r$, the lower bound equals
\begin{equation*}
\frac{T}{2\pi} \cdot \frac{\vartheta r}{1 + \vartheta r} + o(T).
\end{equation*}
Given any $\varepsilon_0>0$, we may first choose $m$ so that $r_m>1-\varepsilon_0$ and then choose $\vartheta<1$ sufficiently close to $1$. It follows that the coefficient of $T$ can be made arbitrarily close to
\[
\frac{1}{2\pi}\cdot\frac{1}{2}=\frac{1}{4\pi}.
\]
Absorbing the remaining difference and the $o(T)$ term into $\varepsilon$ completes the proof of Theorem~\ref{thm:NgCon}.

\section{Higher moments}\label{sect:higher}

We now prove the lower bound of the correct logarithmic order for rational exponents, Theorem~\ref{thm:higher-rational}, together with Corollary~\ref{cor:higher-real-epsilon}. The framework is the one used for the endpoint, and we work on the same range $0<\gamma\le T$. As at the endpoint, the simplicity of the zeros enters only through the residues of $\mathcal R(s)$ at the zeros, that is, through $\mathcal S_1(T)$ below; the mean value $\mathcal S_2(T)$ needs no such hypothesis, since a multiple zero only enlarges it. A mollified first moment is combined with a mollified second moment, and each is evaluated by the residue theorem on a rectangle with vertical sides at $\Re s=c$ and $\Re s=1-c$. The three non-principal edges are treated as before, with the top edge controlled at a good ordinate, the bottom edge of bounded length handled trivially, and the left edge controlled by Lemma~\ref{lem:osc-edge}.

\subsection{The mollified H\"older inequality}\label{ssec:holder-higher}

Fix a rational number $k=a/b\ge1$ in lowest terms, with integers $a\ge b\ge1$. Fix also a parameter $\vartheta$ satisfying
\[
 0<\vartheta<1
 \qquad\text{and}\qquad
 \vartheta\Bigl(1-\frac ba\Bigr)<\frac12.
\]
Since $a\ge b\ge1$ gives $0\le1-b/a<1$, every $\vartheta<\tfrac12$ satisfies both conditions, for every $k$, and the second condition is vacuous when $a=b$. We fix $\vartheta=\tfrac14$ once and for all, so that all implied constants below depend on $k$ alone. Put
\[
  x=T^{\vartheta/a},\qquad X=x^a=T^\vartheta.
\]

Define the short Dirichlet polynomial
\[
  A(s)=\sum_{n\le x}\frac{\lambda(n)d_{1/b}(n)}{n^s}.
\]
This is the unweighted analogue, for the fractional power $\mathcal R(s)^{1/b}$, of the mollifier $\mathcal L_\xi^\dagger(s)$ used in the endpoint case. The coefficients are real, so under the Riemann Hypothesis $A(1-\rho)=\overline{A(\rho)}$.

Choose $T'\in\mathcal T$ with $T-2<T'\le T$, as in Remark~\ref{rem:good-height}, and let $\Gamma_T$ be the set of zeros $\rho$ with $0<\gamma\le T'$. Set
\begin{align}
  \mathcal S_1(T)
  &=\sum_{\rho\in\Gamma_T}
    \frac{\zeta(2\rho)}{\zeta'(\rho)}
    A(\rho)^{a-b}A(1-\rho)^a,\label{eq:S1-higher}\\
  \mathcal S_2(T)
  &=\sum_{\rho\in\Gamma_T}|A(\rho)|^{2a}.
  \label{eq:S2-higher}
\end{align}
Under the Riemann Hypothesis $|A(1-\rho)|=|A(\rho)|$, so each summand of $\mathcal S_1(T)$ has modulus $|\zeta(2\rho)/\zeta'(\rho)|\,|A(\rho)|^{2a-b}$. Applying H\"older's inequality over $\Gamma_T$ with exponents $2a/b=2k$ and $2a/(2a-b)$, and then using the inclusion $\Gamma_T\subseteq\{\rho:0<\gamma\le T\}$ together with the non-negativity of the summands, we obtain
\begin{equation}\label{eq:holder-higher}
  \sum_{0<\gamma\le T}
  \left|\frac{\zeta(2\rho)}{\zeta'(\rho)}\right|^{2k}
  \ge
  \frac{|\mathcal S_1(T)|^{2k}}{\mathcal S_2(T)^{2k-1}},
\end{equation}
valid whenever $\mathcal S_2(T)>0$. This last condition is automatic once Proposition~\ref{prop:S1-higher} is known, since $\mathcal S_1(T)\neq0$ forces $A(\rho)\neq0$ for at least one $\rho\in\Gamma_T$. That proposition is proved in Section~\ref{ssec:S1-higher}, independently of \eqref{eq:holder-higher}.

To prove Theorem~\ref{thm:higher-rational}, we must estimate $\mathcal S_1(T)$ and $\mathcal S_2(T)$. In Sections~\ref{ssec:S1-higher} and~\ref{ssec:S2-higher}, respectively, we establish the following estimates.

\begin{proposition}\label{prop:S1-higher}
Assume the Riemann Hypothesis and that the non-trivial zeros of $\zeta(s)$ are simple. With the notation above,
\[
 |\mathcal S_1(T)|\gg_kT(\log T)^{k^2}.
\]
\end{proposition}

\begin{proposition}\label{prop:S2-higher}
Assume the Riemann Hypothesis and that the non-trivial zeros of $\zeta(s)$ are simple. With the notation above,
\[
 \mathcal S_2(T)\ll_kT(\log T)^{k^2+1}.
\]
\end{proposition}

The exponent in Proposition~\ref{prop:S1-higher} is the true size of $\mathcal S_1(T)$, attained by its diagonal. The bound in Proposition~\ref{prop:S2-higher} is of the order of the term $K_B$ in its proof. In Section~\ref{ssec:higher-completion} we substitute these estimates into \eqref{eq:holder-higher} to prove Theorem~\ref{thm:higher-rational}. The two exponents combine through the identity
\[
 2k\cdot k^2-(2k-1)(k^2+1)=(k-1)^2.
\]

\subsection{Proof of Proposition~\ref{prop:S1-higher}}\label{ssec:S1-higher}

Let $c=1+1/\log T$. The residue theorem applied to
\[
  \mathcal R(s)A(s)^{a-b}A(1-s)^a
\]
over the rectangle with vertical sides $c$ and $1-c$ and horizontal sides $1$ and $T'$ gives exactly $\mathcal S_1(T)$. As in Section~\ref{sect:M1}, the only poles of the integrand inside the rectangle are the non-trivial zeros: there are no zeros with $0<\gamma\le1$, and the pole of $\zeta(2s)$ at $s=\tfrac12$ lies on the real axis, below the bottom edge. Since $[1-c,c]\subset[-1,2]$, the bounds of Lemma~\ref{lem:common-height} apply on the whole of the horizontal edge at the good height $T'$; the bottom edge, at height $1$, is treated by compactness below.

\subsection*{The right edge}
Write
\[
  \mathcal R(s)A(s)^{a-b}
   =\sum_{n=1}^{\infty}\frac{\lambda(n)u(n)}{n^s},
  \qquad
  A(s)^a=\sum_{n\le X}\frac{\lambda(n)v(n)}{n^s},
\]
the common Liouville sign having been extracted by complete multiplicativity. Explicitly,
\[
 u(n)=\sum_{\substack{m\,l_1\cdots l_{a-b}=n\\ l_i\le x}}\ \prod_{i=1}^{a-b}d_{1/b}(l_i),
 \qquad
 v(n)=\sum_{\substack{l_1\cdots l_{a}=n\\ l_i\le x}}\ \prod_{i=1}^{a}d_{1/b}(l_i).
\]
The coefficients $u(n)$ and $v(n)$ are therefore non-negative. They satisfy
\begin{equation}\label{eq:uv-bounds}
  0\le u(n),v(n)\le d_k(n),
  \qquad
  u(n)=v(n)=d_k(n)\quad(n\le x),
\end{equation}
by comparison with the full coefficient identity \eqref{eq:Rz-coeff}, since $d_1*d_{1/b}^{*(a-b)}=d_{1/b}^{*a}=d_k$. Indeed, dropping the truncations only increases the sums, and when $n\le x$ none of the truncations is active.

Write $I_{\mathrm R}$ for the right-edge contribution
\[
 I_{\mathrm R}=\frac1{2\pi i}\int_{c+i}^{c+iT'}\mathcal R(s)A(s)^{a-b}A(1-s)^a\,ds.
\]
Lemma~\ref{lem:MV}, applied at the heights $T'$ and $1$ and differencing, gives
\begin{align}
  I_{\mathrm R}
  &=\frac{T'-1}{2\pi}
    \sum_{n\le X}\frac{u(n)v(n)}n+E_{\mathrm R}.\label{eq:right-higher}
\end{align}
By Lemma~\ref{lem:dk-square} and \eqref{eq:uv-bounds},
\[
  \sum_{n\le X}\frac{u(n)v(n)}n
  \ge\sum_{n\le x}\frac{d_k(n)^2}{n}
  \gg_k(\log T)^{k^2}.
\]
Moreover,
\begin{align*}
 E_{\mathrm R}
 &\ll
 \left(\sum_{n=1}^{\infty}\frac{u(n)^2}{n^{2c-1}}\right)^{1/2}
 \left(\sum_{n\le X}v(n)^2n^{2c-1}\right)^{1/2}
 \ll_k X^c(\log T)^{k^2}
 =o_k\bigl(T(\log T)^{k^2}\bigr),
\end{align*}
where we used Lemma~\ref{lem:dk-square} and $X=T^\vartheta$ with $\vartheta<1$. Since the diagonal in \eqref{eq:right-higher} is real and non-negative, taking real parts gives
\begin{equation}\label{eq:right-higher-main}
  \Re I_{\mathrm R}\gg_kT(\log T)^{k^2}.
\end{equation}

\subsection*{The horizontal edges}
Since $0<1/b\le1$, every value $d_{1/b}(p^j)=\binom{1/b+j-1}{j}=\prod_{i=0}^{j-1}\frac{1/b+i}{1+i}$ lies in $[0,1]$, and so $0\le d_{1/b}(n)\le1$ for all $n$ by multiplicativity. Hence, uniformly for $1-c\le\sigma\le c$ and all $t$,
\begin{equation}\label{eq:A-strip-bound}
 |A(\sigma+it)|
 \le\sum_{n\le x}\frac{d_{1/b}(n)}{n^{\sigma}}
 \le\sum_{n\le x}\frac1{n^{\sigma}}
 \ll x^{\max(1-\sigma,0)}(1+\log x).
\end{equation}
Consequently,
\begin{equation}\label{eq:A-product-horizontal}
 |A(s)^{a-b}A(1-s)^a|
 \ll_k
 x^{(a-b)\max(1-\sigma,0)+a\max(\sigma,0)}(\log T)^{C_k},
\end{equation}
where here and below $C_k$ denotes a positive constant depending only on $k$, whose value may change from line to line.

At the good height $T'$, Lemma~\ref{lem:common-height} and the functional equation give
\begin{equation}\label{eq:R-horizontal-piecewise}
 |\mathcal R(\sigma+iT')|\ll
 \begin{cases}
  T^{1/2+o(1)},&1-c\le\sigma\le0,\\
  T^{1/2-\sigma+o(1)},&0\le\sigma\le\tfrac12,\\
  T^{o(1)},&\tfrac12\le\sigma\le c.
 \end{cases}
\end{equation}
For $\sigma<0$, this follows by reflecting $\zeta(2s)$ at the point $\sigma+iT'$; since $\sigma\ge-1/\log T$, the additional factor $T^{-2\sigma}$ is bounded. For $0\le\sigma\le\tfrac12$ it is the convexity bound for $\zeta(2s)$, and for $\sigma\ge\tfrac12$ it follows from the classical bounds on and to the right of the line $\Re(2s)=1$. The reciprocal factor $1/\zeta(\sigma+iT')$ is $T^{o(1)}$ throughout by \eqref{eq:good-recip}, which is available because $T'$ is a good ordinate and $[1-c,c]\subset[-1,2]$.

We combine \eqref{eq:A-product-horizontal} and \eqref{eq:R-horizontal-piecewise}. If $1-c\le\sigma\le0$, the integrand is
$\ll T^{1/2+\vartheta(1-b/a)+o(1)}$. If $0\le\sigma\le\tfrac12$, its exponent of $T$ is at most
\[
 \frac12+\vartheta\left(1-\frac ba\right)
 +\sigma\left(\frac{\vartheta b}{a}-1\right)+o(1),
\]
which is largest at $\sigma=0$. Finally, for $\tfrac12\le\sigma\le c$ the integrand is
$\ll T^{\vartheta+o(1)}$. The horizontal side at height $T'$ therefore contributes
\begin{equation}\label{eq:horizontal-higher}
 \ll_{k}
 T^{1/2+\vartheta(1-b/a)+o(1)}+T^{\vartheta+o(1)}=o(T),
\end{equation}
since $\vartheta(1-b/a)<\tfrac12$ and $\vartheta<1$.

The bottom edge, of bounded length at height $1$, is handled exactly as the path $I_4$ in Section~\ref{sect:M1}. There $\mathcal R(s)$ is uniformly bounded on the relevant fixed compact subset of the line $\Im s=1$, while \eqref{eq:A-strip-bound} gives $|A(s)^{a-b}A(1-s)^a|\ll_kX(\log T)^{C_k}$. This edge therefore contributes $O_k(X(\log T)^{C_k})=o(T)$.

\subsection*{The left edge}
With $s=1-c+it$, reflect the ratio as in Section~\ref{sect:M1}:
\[
 \mathcal R(s)=\mathcal G(t)\,
 \frac{\zeta(1-2s)}{\zeta(1-s)},
\]
with $\mathcal G(t)$ as in \eqref{eq:Gdef}. Since $\Re(1-2s)=2c-1>1$ and $\Re(1-s)=c>1$, the reflected factors have absolutely convergent expansions
\[
 \zeta(1-2s)=\sum_{q=1}^{\infty}q^{-(2c-1)}(q^2)^{it},
 \qquad
 \frac1{\zeta(1-s)}=\sum_{r=1}^{\infty}\mu(r)r^{-c}r^{it}.
\]
Write also
\[
 A(s)^{a-b}=\sum_{\ell\le x^{a-b}}g_{\ell}\,\ell^{c-1}\ell^{-it},
 \qquad
 A(1-s)^a=\sum_{m\le X}h_m m^{-c}m^{it}.
\]
The coefficients have a common Liouville sign, but only their absolute values are needed here. From the convolution definitions and the non-negativity of $d_{1/b}$,
\begin{align*}
 \sum_{\ell\le x^{a-b}}|g_\ell|\ell^{c-1}
 &\le\left(\sum_{n\le x}d_{1/b}(n)n^{c-1}\right)^{a-b}
 \le\bigl(x^{c}\bigr)^{a-b}
 \ll_k x^{a-b},\\
 \sum_{m\le X}|h_m|m^{-c}
 &\le\left(\sum_{n\le x}\frac{d_{1/b}(n)}{n^c}\right)^a
 \le\zeta(c)^{a/b}
 \ll_k(\log T)^k,
\end{align*}
using $d_{1/b}(n)\le1$ in the first line and $\sum_{n=1}^{\infty}d_{1/b}(n)n^{-c}=\zeta(c)^{1/b}$ in the second.
Therefore the product multiplying $\mathcal G(t)$ is an absolutely convergent sum of the oscillations
\[
 \left(\frac{q^2rm}{\ell}\right)^{it}
\]
over positive rational frequencies. Its total coefficient $\ell^1$-norm is
\begin{align}
 &\left(\sum_{q=1}^{\infty}q^{-(2c-1)}\right)
 \left(\sum_{r=1}^{\infty}\frac{|\mu(r)|}{r^c}\right)
 \left(\sum_\ell|g_\ell|\ell^{c-1}\right)
 \left(\sum_m|h_m|m^{-c}\right)\notag\\
 &\hspace{25mm}\ll_k
 \zeta(2c-1)\frac{\zeta(c)}{\zeta(2c)}
 x^{a-b}(\log T)^{C_k}
 \ll_k x^{a-b}(\log T)^{C_k}.\label{eq:left-l1-higher}
\end{align}

The part of each integral over $[1,t_0]$ is $O(1)$, exactly as in Section~\ref{sect:M1}. Lemma~\ref{lem:osc-edge}, applied with $R=q^2rm/\ell$ over $[t_0,T']$, bounds the rest. Hence
\begin{equation}\label{eq:left-osc-higher}
 \int_{1}^{T'}\mathcal G(t)R^{it}\,dt\ll\sqrt T
 \qquad(R>0),
\end{equation}
uniformly in $R$. Combining \eqref{eq:left-l1-higher} and \eqref{eq:left-osc-higher}, the left edge is
\begin{equation}\label{eq:left-higher}
 \ll_k T^{1/2}x^{a-b}(\log T)^{C_k}
 =T^{1/2+\vartheta(1-b/a)+o(1)}=o(T),
\end{equation}
since $\vartheta(1-b/a)<\tfrac12$.

\subsection*{Conclusion}
Combining \eqref{eq:right-higher-main}, \eqref{eq:horizontal-higher}, and \eqref{eq:left-higher}, and taking real parts, gives
\[
 |\mathcal S_1(T)|\ge \Re\mathcal S_1(T)
 =\frac{T'-1}{2\pi}\sum_{n\le X}\frac{u(n)v(n)}n
 +o_k\!\left(T(\log T)^{k^2}\right)
 \gg_kT(\log T)^{k^2}.
\]
This proves Proposition~\ref{prop:S1-higher}.

\subsection{Proof of Proposition~\ref{prop:S2-higher}}\label{ssec:S2-higher}

Put
\[
 B(s)=A(s)^a=\sum_{n\le X}\frac{b_n}{n^s},
 \qquad b_n=\lambda(n)v(n),
 \qquad 0\le v(n)\le d_k(n),
\]
and
\[
 D_B=\sum_{n\le X}\frac{v(n)^2}{n},
\]
so that Lemma~\ref{lem:dk-square} gives
\begin{equation}\label{eq:DB-bound}
 D_B\ll_k(\log T)^{k^2}.
\end{equation}
The coefficients $b_n$ are real, so on the Riemann Hypothesis $B(1-\rho)=\overline{B(\rho)}$ and $|A(\rho)|^{2a}=B(\rho)B(1-\rho)$. Under simplicity the residue theorem returns exactly one term per zero. Hence
\[
 \mathcal S_2(T)=\frac1{2\pi i}\oint
 \frac{\zeta'(s)}{\zeta(s)}B(s)B(1-s)\,ds
\]
over the same rectangle.

On either horizontal side, partial summation and the usual summatory estimate for $d_k$ (see \cite[Chapter~II.5]{Tenenbaum}) give
\[
 |B(\sigma+it)|\ll_k
 X^{\max(1-\sigma,0)}(\log T)^{C_k},
\]
and hence
\[
 |B(\sigma+it)B(1-\sigma-it)|
 \ll_kX(\log T)^{C_k}
 \qquad(1-c\le\sigma\le c).
\]
At the height $T'$ we use \eqref{eq:good-logder}, which covers the whole edge since $[1-c,c]\subset[-1,2]$. On the bottom path, of bounded length at height $1$, the factor $\zeta'(s)/\zeta(s)$ is uniformly bounded, as in Section~\ref{sect:M2}. The two horizontal contributions are therefore
\[
 O_k\bigl(X(\log T)^{C_k}\bigr)=o(T).
\]

On the left edge we use
\[
 \frac{\zeta'(1-s)}{\zeta(1-s)}
 =-\frac{\zeta'(s)}{\zeta(s)}-\log\frac{t}{2\pi}+O(t^{-1}).
\]
Change variables $s\mapsto1-s$ and account for the downward orientation. Taking complex conjugates is then legitimate, because the $b_n$ and the Dirichlet coefficients of $\zeta$ are real, so that $B(\bar s)=\overline{B(s)}$ and $\zeta'(\bar s)/\zeta(\bar s)=\overline{\zeta'(s)/\zeta(s)}$. The right edge is exactly $J_B'$, and adding the two gives
\begin{equation}\label{eq:S2-KJ}
 \mathcal S_2(T)=\overline{K_B}+2\Re J_B'
 +O_k\bigl(X(\log T)^{C_k}\bigr),
\end{equation}
where
\[
\begin{aligned}
 K_B&=\frac1{2\pi}\int_{1}^{T'}\log\frac{t}{2\pi}\,B(c+it)B(1-c-it)\,dt,\\
 J_B'&=\frac1{2\pi}\int_{1}^{T'}\frac{\zeta'(c+it)}{\zeta(c+it)}B(c+it)B(1-c-it)\,dt.
\end{aligned}
\]
The displayed error includes the horizontal sides and the $O(t^{-1})$ term. The latter is bounded using the trivial estimate
$|B(c+it)B(1-c-it)|\ll_kX(\log T)^{C_k}$.

\subsection*{The term $K_B$}
For $1\le u\le T'$, set
\[
 H_B(u)=\int_{1}^{u}B(c+it)B(1-c-it)\,dt.
\]
Lemma~\ref{lem:MV}, applied at the heights $u$ and $1$ and differencing, gives, uniformly in $u$,
\begin{equation}\label{eq:HB-MV}
 H_B(u)=(u-1)D_B+O_k\bigl(X^c(\log T)^{k^2}\bigr).
\end{equation}
Indeed,
\[
\begin{aligned}
 \sum_{n\le X}v(n)^2n^{1-2c}
 &\le\sum_{n=1}^{\infty}\frac{d_k(n)^2}{n^{1+2(c-1)}}
 \ll_k(\log T)^{k^2},\\
 \sum_{n\le X}v(n)^2n^{2c-1}
 &\le X^{2c}\sum_{n\le X}\frac{v(n)^2}{n}
 \ll_kX^{2c}(\log T)^{k^2}.
\end{aligned}
\]
Partial summation in $t$, using \eqref{eq:HB-MV} and \eqref{eq:DB-bound}, yields
\begin{equation}
 K_B
 \ll T(\log T)D_B+X(\log T)^{C_k}
 \ll_kT(\log T)^{k^2+1}.
 \label{eq:KB-bound}
\end{equation}

\subsection*{The term $J_B'$}
Write, in the notation of Section~\ref{sect:M2} but now with $B$ in place of $\mathcal L^\dagger_\xi$,
\[
 -\frac{\zeta'(s)}{\zeta(s)}B(s)
 =\sum_{n=1}^{\infty}\frac{\alpha_n}{n^s},
 \qquad
 \alpha_n=\sum_{\substack{hm=n\\ m\le X}}\Lambda(h)b_m.
\]
Since $|b_m|=v(m)\le d_k(m)$,
\begin{equation}\label{eq:alpha-general-bound}
 |\alpha_n|
 \le\sum_{h\mid n}\Lambda(h)\,d_k(n/h)
 =\frac{d_k(n)\log n}{k};
\end{equation}
the identity follows by comparing coefficients after logarithmically differentiating $\zeta(s)^k$. Lemma~\ref{lem:MV}, applied at the heights $T'$ and $1$ and differencing, now gives
\begin{equation}\label{eq:JB-MV}
 J_B'=-\frac{T'-1}{2\pi}
 \sum_{n\le X}\frac{\alpha_nb_n}{n}+E_B.
\end{equation}
By \eqref{eq:alpha-general-bound} and Lemma~\ref{lem:dk-square},
\begin{equation}\label{eq:JB-diagonal}
 \sum_{n\le X}\frac{|\alpha_nb_n|}{n}
 \ll_k\sum_{n\le X}\frac{d_k(n)^2\log n}{n}
 \ll_k(\log T)^{k^2+1}.
\end{equation}
For the error term,
\begin{align*}
 E_B
 &\ll
 \left(\sum_{n=1}^{\infty}\frac{|\alpha_n|^2}{n^{2c-1}}\right)^{1/2}
 \left(\sum_{n\le X}v(n)^2n^{2c-1}\right)^{1/2}\\
 &\ll_k
 \bigl((\log T)^{k^2+2}\bigr)^{1/2}
 \bigl(X^{2c}(\log T)^{k^2}\bigr)^{1/2}
 \ll_k X(\log T)^{k^2+1}=o_k(T),
\end{align*}
where the first factor follows from \eqref{eq:alpha-general-bound} and the $j=2$ case of Lemma~\ref{lem:dk-square}. Equations \eqref{eq:JB-MV} and \eqref{eq:JB-diagonal} therefore imply
\begin{equation}\label{eq:JB-bound}
 J_B'\ll_kT(\log T)^{k^2+1}.
\end{equation}
Combining \eqref{eq:S2-KJ}, \eqref{eq:KB-bound}, and \eqref{eq:JB-bound} proves Proposition~\ref{prop:S2-higher}.

\subsection{Completion of the proofs}\label{ssec:higher-completion}

Inserting Propositions~\ref{prop:S1-higher} and~\ref{prop:S2-higher} into \eqref{eq:holder-higher} gives
\begin{align*}
 \sum_{0<\gamma\le T}
 \left|\frac{\zeta(2\rho)}{\zeta'(\rho)}\right|^{2k}
 \gg_k
 \frac{\{T(\log T)^{k^2}\}^{2k}}
      {\{T(\log T)^{k^2+1}\}^{2k-1}}
 =T(\log T)^{(k-1)^2}.
\end{align*}
This proves Theorem~\ref{thm:higher-rational}.

\begin{proof}[Proof of Corollary~\ref{cor:higher-real-epsilon}]
The case $k=1$ follows from Theorem~\ref{thm:NgCon}. Let $k>1$, write $w_\rho=\zeta(2\rho)/\zeta'(\rho)$, and choose a rational $r$ with $1\le r<k$ and $k(k-r)<\varepsilon$. H\"older's inequality and the bound $N(T)\ll T\log T$ give
\[
 \sum_{0<\gamma\le T}|w_\rho|^{2r}
 \le
 \left(\sum_{0<\gamma\le T}|w_\rho|^{2k}\right)^{r/k}
 N(T)^{1-r/k}.
\]
Theorem~\ref{thm:higher-rational} at the exponent $r$, followed by rearrangement, yields
\begin{align*}
 \sum_{0<\gamma\le T}|w_\rho|^{2k}
 \gg_{k,r}T(\log T)^{kr-2k+1}
 =T(\log T)^{(k-1)^2-k(k-r)}.
\end{align*}
Since $k(k-r)<\varepsilon$ and $\log T\ge1$, the exponent exceeds $(k-1)^2-\varepsilon$, which is the claim.
\end{proof}

\begin{remark}\label{rem:rationality}
The restriction to rational $k$ in Theorem~\ref{thm:higher-rational} enters only through the representation $k=a/b$ and the use of integer powers of the polynomial $A(s)$. Removing it without the $\varepsilon$-loss should be possible by adapting the truncated-exponential machinery used by Gao and Zhao for real negative moments of $\zeta'(\rho)$ \cite{GaoZhao}. We do not carry out that separate argument here.
\end{remark}

\section*{Acknowledgements}
Thanks to Ben Durkan, Chris Hughes and Nathan Ng for carefully reading earlier drafts of this paper and for their helpful suggestions. Thanks also to Saloni Sinha for the correspondence regarding an earlier version of this paper. The author gratefully acknowledges support from the Heilbronn Institute for Mathematical Research.

%
%
%

\end{document}